\documentclass[11pt,
]
{article}

\usepackage{amssymb,amstext,amsthm,latexsym}
\usepackage
{amsmath}
\usepackage{amsfonts}
\usepackage{bbm}

\usepackage{mathrsfs}
\usepackage{mathbbol}

\usepackage{mparhack}

\input{im.sty}

\begin{document}

\title
{A generic property of Solovay's set $\Sigma$ 
}

\author 
{
Vladimir~Kanovei\thanks{IITP RAS and MIIT,
  Moscow, Russia, \ {\tt kanovei@googlemail.com} --- contact author. 
}  
\and
Vassily~Lyubetsky\thanks{IITP RAS,
  Moscow, Russia, \ {\tt lyubetsk@iitp.ru} 
}
}


\date{\today}
\maketitle


\begin{abstract}
We prove that Solovay's set $\Sg$ is generic over the 
ground model via a forcing notion whose order relation 
\ddd\sq extends the given order relation.%
\end{abstract}

\def\contentsname{}


\punk{Introduction}
\las{x1}

Solovay's paper \rit{A model of set theory in which every 
set of reals is Lebesgue measurable} \cite{sol} 
belongs to the 
classics of forcing. 
Its main reasult is the construction of the said model,  
and a related model in which only definable sets are 
claimed to be Lebesgue measurable, but the axiom of choice 
holds, unlike the titular model, where it fails by 
necessity. 
During the course of the paper, Solovay invented several 
cornerstone methods in the forcing practice. 
One of those is an ``important lemma'' of 
\cite[Section 4.4]{sol}, which asserts, roughly, that a 
generic extension of a ground model also is a generic 
extension of any intermediate submodel via a forcing 
notion $\Sg$ which is a subset of the original 
forcing notion $\dP$. 

We prove in this note (Theorem~\ref{sig2}), 
that, in the context of the Solovay construction, the set 
$\Sg$ itself is generic over the ground model, 
via a forcing notion closely related to $\dP$ in the 
sense that its domain is equal to the domain of the 
original forcing $\dP$ but the partial order relation 
\ddd\sq extends the given partial order relation 
$\le_\dP$.

As an application, we present a short proof of 
the following result (Theorem~\ref{is}): 
any subextension of a 
Cohen-generic extension is a Cohen-generic extension 
itself (or else trivial), and the whole extension 
is a Cohen-generic extension 
of the subextension (or else trivial). 
This is a folklore result known  
probably from the first years 
of forcing, but in fact it does not seem to have 
ever been published with proof.

\punk{The intermediate forcing $\Sg$}
\las{x2}

The exact content of Solovay's result is as follows:

\bpro
\lam{st}
Suppose that\/ $\dP\in\vv$ is a forcing notion in the 
set universe\/ $\vv$, $G\sq\dP$ a\/ \ddd\dP generic 
filter over\/ $\vv$, $t\in\vv$ a\/ \ddd\dP name, 
$X=t[G]\sq\vv$\snos
{$t[G]$ is the \ddd G\rit{valuation}, 
also called \rit{interpretation}, of the name $t$; 
$t[G]\in\vv[G]$.}. 
Then there is a set\/ 
$\Sg=\xha Xt\in\vv[G]$, $\Sg\sq\dP$ such that
\ben
\Renu
\vyk{
\itla{s1}
the inequality 
$\xha Xt\ne\pu$ is necessary and sufficient for there to exist 
a set $G\sq\dP$, \ddd\dP generic over $\vv$ and satisfying 
$X=t[G]$;
}%
\itla{st1}
$\Sg$ is closed weakwards in\/ $\dP$, so that 
if\/ $q\in\Sg$, $p\in\dP$, and\/ $p\le q$, 
then\/ $q\in\Sg\,;$\snos
{We always assume that $p\le q$ means that $q$ is a 
stronger condition.}

\itla{st2}
$\vv[\Sg]=\vv [X]\,;$ 

\itla{st3}
$G\sq\Sg$ and\/ $G$ is\/ 
\ddd{\Sg}generic over\/ $\vv[X]\,;$ 

\itla{st4}
therefore\/ $\vv[G]$ is a\/ \ddd{\Sg}generic 
extension of\/ $\vv[X]=\vv[\Sg]\,;$ 

\itla{st5}
if a set\/ $G'\sq\Sg$ is \ddd{\Sg}generic 
\imar{st5}
over $\vv[X]$ then\/ $G'$ is\/ \ddd\dP generic 
over\/ $\vv$ and still\/ $t[G']=X$.
\een
\epro
\bpf
[sketch, see detailed arguments in \cite{sol}, 4.4]
Let $\Sg=\dP\bez\bigcup_{\xi<\vt}A_\xi$, where the 
sequence of sets $A_\xi\sq\dP$ is defined in $\vv[G]$ 
as follows:
\ben
\aenu
\itla{si1}
$A_0$ 
\imar{si1}
consists of all conditions $p\in\dP$ which either force 
$\dtx \in t$ for some $x\in\vv\bez X$, or force 
$\dtx \nin t$ for some $x\in X$.\snos
{If $x\in\vv$ then $\dtx$ is a canonical name for $x$.}

\itla{si2}
$A_{\xi+1}$ 
consists of all conditions $p\in\dP$ such that there is 
a dense set $D\in\vv$ in $\dP$ satisfying: 
if $q\in D$ and $p\le q$ then $q\in A_\xi$.

\itla{si3}
$A_\la=\bigcup_{\xi<\la}A_\xi$ whenever $\la$ is a limit 
ordinal.
\een
To conclude, each condition $p\in A_0$  
directly contradicts the assumption 
that $t$ is a name for $X$, by \ref{si1}, and 
this contradiction prevails by \ref{si2} and 
\ref{si3} for all bigger $\xi$ in more and 
more indirect way.
This \ddd\sq increasing sequence of sets $A_\xi\sq\dP$ 
stabilizes on a limit ordinal $\vt\in\vv$, and we have 
got the sets 
$A=\bigcup_{\xi<\vt}A_\xi$ and $\Sg=\dP\bez A$. 
\epf

Further studies of Grigorieff \cite{gri} and others 
on intermediate submodels of generic extensions 
demonstrated that not only $\vv[G]$ is a generic extension 
of $\vv[X]$ by \ref{st3}, but $\vv[X]$ itself is a 
generic extension of $\vv$.\snos 
{See also \cite[Fact 11]{fhr}, \cite[15.43]{jech}, 
\cite[Proposition 10.10]{kanam}, \cite{kl21}, 
or \cite[Section 1]{zapt}, among other references.}
Theorem~\ref{sig2} below, our main result, 
asserts that \rit{the set\/ 
$\xha Xt$ itself is a generic filter over\/ $\vv$}, 
via a forcing notion closely related to $\dP$.

\punk{The genericity of the set $\Sg$}
\las{x3}

Arguing in the context of Proposition~\ref{st}, 
we define a new order $\le_t$ on $\dP$, which extends the 
original order ${\le}={\le_\dP}$, as follows: 
$p\le_t q$ iff $q$ $\dP$-forces over $\vv$ that 
$\dtp\in \xha{\dtt[\doG]}{\dtt}$. 
In other words, for $p\le_t q$ it is necessary and sufficient 
that we have $p\in\xha Xt$ whenever $G\sq\dP$ is generic over 
$\vv$, $X=t[G]$, and $q\in G$.

\ble
\lam{sig1}
$\le_t$ is a partial (pre)order relation on\/ $\dP$ 
which belongs to\/ $\vv$ and extends the given 
order\/ ${\le}={\le_\dP}$, 
so that\/ ${\le_\dP}\sq {\le_t}$ 
(or equivalently, $p\le_\dP q$ implies\/ $p\le_t q$).
\ele
\bpf
Suppose that $p\le_t q\le_t r$. 
To prove $p\le_t r$, assume that $G\sq\dP$ is generic 
over $\vv$, $r\in G$, and $X=t[G]$; 
we have to prove that $p\in \xha Xt$.

By definition, $q\in \xha Xt$.
Pick a set $G'\sq  \xha Xt$ \ddd{\xha Xt}generic over $\vv[X]$ 
and containing $q$.
Then $G'$ is \ddd\dP generic over $\vv$ and still $t[G']=X$, 
by \ref{st5}.
Therefore $p\in \xha Xt$, because $p\le_t q$, and we are done.

Finally, suppose that $p\le q$ and prove $p\le_t q$. 
Assume that $G\sq\dP$ is generic 
over $\vv$, $q\in G$, and $X=t[G]$. 
To prove that $p\in \xha Xt$, note first of all that 
$q\in \xha Xt$ by \ref{st3}; 
then $p\in \xha Xt$ by \ref{st1}.
\epf

\bte
[in the assumptions of Proposition~\ref{st}]
\lam{sig2}
Suppose that\/ $G\sq\dP$ is generic over\/ $\vv$ 
and\/ $X=t[G]$. 
Then the set\/ $\Sg=\xha Xt$ itself is a generic filter 
over\/ $\vv$ 
in the forcing\/ $\dP_t=\stk{\dP}{\le_t}$.
\ete
\bpf
Suppose that $p\in \dP$, $q\in \Sg$, and $p\le_t q$. 
To prove $p\in \Sg$, 
consider any set $G'\sq  \Sg$, 
\ddd{\Sg}generic over $\vv[X]$ and containing $q$. 
Then $G'$ is \ddd\dP generic over $\vv$, and $t[G']=X$, 
by \ref{st5}.
Now we have $p\in \Sg$ since $p\le_t q$.

Prove that any two conditions $p,q\in \Sg$ are 
\ddd{\le_t}compatible in the set $\Sg$. 
By genericity there is a condition $r\in G$ which forces 
both $\dtp\in \xha{\dtt[\doG]}{\dtt}$ and 
$\dtq\in \xha{\dtt[\doG]}{\dtt}$. 
Then by definition $p\le_t r$ and $q\le_t r$, and 
on the other hand $r\in \Sg$ by \ref{st3}.

Finally prove the genericity itself. 
Suppose that a set $D\sq\dP$ is \ddd{\le_t}dense 
(not necessarily dense \poo\ the original order). 
Assume towards the contrary that some $p\in G$ 
forces $\breve D\cap \xha{\dtt[\doG]}{\dtt}=\pu$. 
By density, there is a condition $q\in D$, $p\le_t q$. 

Consider a set $G'\sq\dP$, \ddd\dP generic over $\vv$ 
and containing $q$, and let $X'=t[G']$. 
Then $q\in G'\sq\Sg'=\xha{X'}t$, and hence $p\in\Sg'$ 
by the above. 

Further, consider a set $G''\sq\Sg'$, 
\ddd{\Sg'}generic over $\vv[X']$ 
and containing $p$. 
Then $G''$ is \ddd\dP generic over $\vv$, and $t[G'']=X'$, 
by \ref{st5}. 
Thus $q\in D\cap\xha{t[G'']}t$ and $p\in G''$,  
contrary to the choice of $p$. 
\epf

\punk{Intermediate submodels of Cohen-generic 
extensions}
\las{x4}

Now we are getting the following for free.

\bte
[forcing folklore]
\lam{is}
Assume that\/ $a\in\dn$ is a Cohen-generic real over 
the set universe\/ $\vv$, and\/ $X\in\vv[a],\:X\sq\vv$. 
Then
\ben
\renu
\itla{is1}
either\/ $X\in\vv$ or\/ $\vv[X]$ is a Cohen-generic 
extension of\/ $\vv\,;$

\itla{is2}
either\/ $\vv[X]=\vv[a]$ or\/ $\vv[a]$ is a Cohen-generic 
extension of\/ $\vv[X]\,.$
\een
\ete
\bpf
Note that the Cohen forcing is countable.
Therefore both $\Sg$ in Proposition~\ref{st} and $\dQ$ in 
Theorem~\ref{sig2} are countable forcing notions.
\epf

\end{document}